\documentclass[11pt]{article}
\usepackage{graphicx} 

\RequirePackage{amsthm,amsmath,amsfonts,amssymb}

\usepackage[
  top=2.5cm,
  bottom=2.5cm,
  left=2.5cm,
  right=2.5cm
]{geometry}

\usepackage{url}
\usepackage{adjustbox}
\include{macros1}
\usepackage{moresize}
\usepackage{scalefnt}
\usepackage{multirow}
\usepackage{amsmath}
\allowdisplaybreaks
\usepackage[font=small]{caption}
\usepackage{xr}
\usepackage{mathrsfs}
\usepackage{natbib}
\usepackage{comment}
\usepackage{array}
\usepackage{xcolor}
\usepackage{graphicx}
\usepackage{algorithm}
\usepackage{algpseudocode}
\usepackage{dsfont}

\def\zero{{\mathbf {0}}}

\def\trans{^{\mbox{\tiny{\sf T}}}}
\newcommand{\real}{\mathbb {R}}

\newcommand\Ldx[1]{\mathcal {B}^{#1}(X)}
\newcommand\Ldfx[1]{\mathcal {B}^{#1}\{f_0(X)\}}
\def\Ls{\mathcal {L}}

\def\var{{\mathrm {Var}}}

\def\inv{^{-1}}

\def\bop{O_P}

\def\var{{\mathrm {var}}}
\def\tr{{\mathrm {tr}}}

\def\K{\mathbb {K}}
\def\U{{\mathbb {U}}}
\def\V{{\mathbb {V}}}
\def\vecone{{\mathds {1}}}
\def\Q{{\mathbb {Q}}}
\def\I{{\mathbb {I}}}
\def\H{{\mathcal {H}}}
\def\G{{\mathbb {G}}}
\def\D{{\mathbb {D}}}
\def\W{{\mathbb {W}}}

\def\Gs{{\mathbf {G}}}
\def\wgsir{w_{{\mathrm {GSIR}}}}

\renewcommand{\indep}{\;\, \rule[0em]{.03em}{.67em} \hspace{-.25em}}

\theoremstyle{thmstyleone}%
\newtheorem{theorem}{Theorem}

\newtheorem{example}{Example}
\theoremstyle{remark}

\graphicspath{{Fig/}}

\title{Supervised Laplacian embedding: bridging dimensionality reduction and sufficient dimension reduction}
\author{Wei Luo, Center for Data Science, Zhejiang University}
\date{November 2026}

\begin{document}

\maketitle

\abstract{Dimensionality reduction in machine learning research has largely developed in isolation from sufficient dimension reduction in the statistical literature. In this paper, we first develop a new supervised Laplacian embedding method that addresses issues such as inappropriate normalization in existing approaches. We then develop statistical theory for the effectiveness of the proposed method, based on a novel, rigorous definition of an appropriate similarity measure that is also central to the framework of Laplacian embedding. Finally, we uncover the intrinsic connection between the proposed method and popular sufficient dimension reduction methods, thereby bridging the gap between machine learning research on dimensionality reduction and statistical research on sufficient dimension reduction for the first time in the literature. This connection benefits both communities in multiple ways: by bringing theoretical rigour to dimensionality reduction, by providing an algorithmic perspective on sufficient dimension reduction, and by generating new methods drawing on each other's literature. An example of the latter, called second-order supervised Laplacian embedding, is discussed in detail. Numerical studies are presented at the end. }


\section{Introduction}
\label{sec: intro}

Dimensionality reduction, often as part of feature engineering or representation learning, is an important data pre-processing technique that de-noises high-dimensional data to obtain a lower-dimensional representation and thereby facilitates subsequent data analysis. Laplacian embedding \citep{belkin2003laplacian,li2019laplacian} is a popular unsupervised dimensionality reduction method that constructs the reduced data via a nonlinear transformation. The essence of Laplacian embedding lies in a measure of similarity between sample observations. Let $X$ be the random object that generates the data $X_1,\ldots, X_n$. The similarity between $X_i$ and $X_j$ is a nonnegative scalar $w_{X, ij}$ that is larger when $X_i$ and $X_j$ are closer under a certain metric, such as the Euclidean distance if they are real vectors. The reduced data from Laplacian embedding, denoted by $U_1,\ldots, U_n \in \real^d$, are defined to best preserve this similarity measure or equivalently to minimise
\begin{equation}\label{eq: le original}
\sum_{i=1}^n \sum_{j=1}^n \|U_i - U_j\|^2 w_{X,ij}
\end{equation}
subject to $\sum_{i=1}^n U_i U_i\trans = n \I_d$, where $\|\cdot\|$ denotes the $L_2$ norm of a real vector and $\I_d$ denotes the $d$-dimensional identity matrix. Let $\W_X$ be the $n\times n$ matrix whose $(i,j)$th entry is $w_{X,ij}$, let $\D_X$ be the $n\times n$ diagonal matrix whose $i$th diagonal entry is $\sum_{j=1}^n w_{X,ij}$, and let $\U = (U_1,\ldots, U_n)\trans$. With $\tr(\cdot)$ denoting the trace of a square matrix, the minimisation in (\ref{eq: le original}) can be written as
\begin{equation}\label{eq: le}
\min_{\U \in \real^{n\times d}} \tr \{\U\trans (\D_X - \W_X) \U \} \quad \mbox{subject to } \U\trans \U = n \I_d.
\end{equation}
Usually, $\D_X - \W_X$ is called the Laplacian matrix. Since we can assume the symmetry of $\W_X$ in (\ref{eq: le}) without loss of generality, (\ref{eq: le}) reduces to a standard eigen-decomposition problem. For this reason, Laplacian embedding is also commonly referred to as Laplacian eigenmap.

In supervised learning, where a response variable $Y$ is of interest, the purpose of dimensionality reduction is to find a lower-dimensional transformation of the predictor $X$ that preserves the information for modelling $Y$. Motivated by the commonly observed effectiveness of (unsupervised) Laplacian embedding, several supervised Laplacian embedding methods have been proposed; see, for example, \cite{costa2005}, \cite{jiang2009}, \cite{raducanu2012}, \cite{gong2014}, and \cite{tai2019}. Naturally, these methods aim to find a reduced predictor whose similarity measure agrees maximally with the similarity measure of $Y$. As with $w_{X,ij}$ above, the latter is a nonnegative scalar $w_{Y, ij}$ that is larger when the corresponding observations $Y_i$ and $Y_j$ are closer. Consequently, the reduced predictor $\U = (U_1,\ldots, U_n)\trans \in \real^{n\times d}$ minimises
\begin{equation}\label{eq: le sup}
 \sum_{i=1}^n \sum_{j=1}^ n \|U_i - U_j\|^2 w_{Y,ij} \quad \mbox{or equivalently} \quad 
   \tr\{\U\trans (\D_Y - \W_Y) \U\},
\end{equation}
where $\W_Y$ and $\D_Y$ are defined analogously to $\W_X$ and $\D_X$ above, with $w_{Y,ij}$ in place of $w_{X,ij}$. As with (\ref{eq: le}), the minimiser of (\ref{eq: le sup}) exists under some normalisation constraint on $U_i$'s and an additional regularisation used to avoid overfitting \citep{tai2019}. 

As pointed out by \cite{tai2022kernelized}, all previous Laplacian embedding methods, both unsupervised and supervised, suffer from an ``out-of-sample problem''. That is, their output consists of observations $U_1,\ldots, U_n$ rather than an explicit transformation of $X$, and therefore cannot directly produce $U_{n+1}$ for a new observation $X_{n+1}$. To address this issue, \cite{tai2022kernelized} proposed a kernelised supervised Laplacian embedding method, assuming that $U_1,\ldots, U_n$ in (\ref{eq: le sup}) are observations of $f(X)$ for some set of functions $f(\cdot)$ belonging to a prespecified reproducing kernel Hilbert space (RKHS). Let $\H_X$ and $K_X(\cdot,\cdot)$ denote the RKHS and the corresponding kernel function, respectively, and let $\K_X \in \real^{n\times n}$ be the matrix whose $(i,j)$th entry is $K_X(X_i,X_j)$. The assumption $f(\cdot) \subset \H_X$ implies $\U = \K_X \V$ for some $\V\in \real^{n\times d}$, which turns (\ref{eq: le sup}) into 
\begin{equation}\label{eq: le sup kernel tai22}
\tr\{\V\trans \K_X (\D_Y - \W_Y) \K_X \V\}.
\end{equation}
To normalise $\V$ in (\ref{eq: le sup kernel tai22}), \cite{tai2022kernelized} adopted the constraint
\begin{equation}\label{eq: le sup kernel tai22 constr}
\V\trans \K_X \D_Y \K_X \V = n \I_d,
\end{equation}
under which the minimisation of (\ref{eq: le sup kernel tai22}) is a generalised eigen-decomposition problem. Let $\ddot \V$ be the corresponding minimiser. The reduced predictor is $(K_X(X,X_1),\ldots, K_X(X,X_n))\ddot \V$. 

Despite its success in resolving the ``out-of-sample problem'', the supervised Laplacian embedding method proposed in \cite{tai2022kernelized} can be ineffective, as the following example illustrates. The reason is that the normalisation constraint (\ref{eq: le sup kernel tai22 constr}) involves $Y$ and is therefore incompatible with the purpose of dimensionality reduction, which is to maximise, rather than normalise, the alignment between the reduced predictor and the response. Originally, \cite{tai2022kernelized} incorporated ``feature information'' into dimensionality reduction by modifying $\D_Y - \W_Y$ in (\ref{eq: le sup kernel tai22}) to $\lambda (\D_X - \W_X) + (1-\lambda) (\D_Y - \W_Y)$ and correspondingly modifying $\D_Y$ in (\ref{eq: le sup kernel tai22 constr}) to $\lambda \D_X + (1-\lambda) \D_Y$, where the prespecified constant $\lambda$ ``controls the extent to which the feature information is prioritized against the label information''. However, this modification does not resolve the ineffectiveness caused by inappropriate normalisation. To the best of our knowledge, there is currently no supervised Laplacian embedding method in the literature that delivers an effective reduced predictor while avoiding the ``out-of-sample problem''. 


\def\ind{{\mathcal {I}}}

\begin{example}\label{example: tai inconsist}
Suppose $Y$ is categorical and $w_{Y,ij} = \ind (Y_i = Y_j)$, where $\ind(\cdot)$ denotes the indicator function. Then the $i$th diagonal element of $\D_Y$ is the sample proportion of the category to which $Y_i$ belongs, which means that the minimisation of (\ref{eq: le sup kernel tai22}) subject to (\ref{eq: le sup kernel tai22 constr}) asymptotically approximates
\begin{equation}\label{eq: tai example}
\min_{f(\cdot) \subset \H_X} \tr\left[E[\var\{f(X)\mid Y\} \pi(Y)] \right] \quad \mbox{ subject to } E\{f(X) f\trans(X) \pi(Y)\} = \I_d,
\end{equation}
where $\pi(\cdot)$ denotes the probability mass function of $Y$.
For simplicity, we focus on linear dimensionality reduction; that is, $\H_X$ consists of linear functions. Then the reduced predictor derived from (\ref{eq: tai example}) is 
$(A \beta)\trans X$, where $A = [E\{X X\trans \pi(Y)\}]^{-1/2}$ and $\beta$ is the set of leading eigenvectors of $A E\{E(X\mid Y) E(X\trans\mid Y) \pi(Y)\} A$. This differs from the desired linear reduced predictor---for example, $B\trans X$ in linear discriminant analysis, where $B$ is the set of leading left singular vectors of $[E\{\var (X \mid Y)\}]^{-1} E\{E(X\mid Y) E(X\trans\mid Y) \}$. Specifically, if $Y$ is binary with $\pi(1) = \pi(2) = 0.5$ and $X\mid Y \sim N((3Y,3)\trans, \I_2)$, then (\ref{eq: tai example}) yields $(0.52, 0.85)\trans X$, whereas the desired reduced predictor is the first variable in $X$.
\end{example}

In this paper, we propose a new supervised Laplacian embedding method, in which we follow \cite{tai2022kernelized} to avoid the ``out-of-sample problem'' by kernelising the reduced predictor, but impose an appropriate normalisation constraint on $\V$ by slightly modifying a given similarity measure $w_{Y,ij}$. 
In addition, we develop statistical theory to justify the effectiveness of the proposed method; that is, it delivers the desired reduced predictor which preserves the information in $X$ for modelling $Y$. As the foundation of this theory, we provide, for the first time in the literature, a rigorous definition of an appropriate similarity measure $w_{Y,ij}$. Surprisingly, this definition does not merely restrict $w_{Y,ij}$ but also relaxes it by allowing $w_{Y,ii} < w_{Y,ij}$, which contradicts the intuition that $w_{Y,ij}$ must increase as $Y_j$ approaches $Y_i$. 

In the statistical literature, there is a body of research on supervised dimensionality reduction, commonly known as sufficient dimension reduction \citep[SDR;][]{li1989,sdr18} or, more specifically, nonlinear SDR \citep{lee13} when a nonlinear transformation is used to generate the reduced predictor. Instead of preserving the similarity of $Y$, existing nonlinear SDR methods focus on characterising the sufficiency of the reduced predictor through a collection of appropriate conditional moments. 
Interestingly, our theory is intrinsically connected with the statistical theory of nonlinear SDR: the proposed definition of a similarity measure is closely related to both the kernel function of an RKHS and a concept known as a continuous negative definite function, and the proposed supervised Laplacian embedding method can be regarded as a generalisation of some state-of-the-art nonlinear SDR methods \citep{lee13, chen2024}. This is the first attempt in the literature to bridge dimensionality reduction in machine learning research and SDR in statistical research, both of which have been extensively developed over decades but are largely separate from each other. By providing a novel statistical formulation of dimensionality reduction and a novel algorithmic perspective on nonlinear SDR methods, this connection will benefit both communities. For example, it can bring theoretical rigour and new insights to dimensionality reduction, provide an alternative and arguably simpler interpretation of nonlinear SDR, and generate new methods by drawing on each other's literature. As an example of the latter, we briefly introduce second-order supervised Laplacian embedding. 


The rest of the article is organised as follows. In Section \ref{sec: sle}, we modify a given similarity measure $w_{Y,ij}$ and propose the corresponding supervised Laplacian embedding method. In Section \ref{sec: sle theory}, we rigorously define an appropriate similarity measure of $Y$ and subsequently develop statistical theory for the effectiveness of the proposed method. Building on these, in Section \ref{sec: sdr} we discuss the connection between the proposed method and several popular nonlinear SDR methods, as well as how this sheds light on the development of both dimensionality reduction and nonlinear SDR. Second-order supervised Laplacian embedding is briefly introduced in Section \ref{sec: second-order sle}. Synthetic and real data analyses are provided in Section \ref{sec: simulation} and Section \ref{sec: real data}, respectively. For continuity of presentation, we defer all theoretical proofs to the Supplementary Material. \texttt{R} code for implementing the proposed method can be downloaded at \url{https://github.com/wluostat/SLE}.

For consistency of presentation, we use statistical language throughout the article, referring to $Y$ as the response rather than the label, and to $X$ as the predictor rather than the features. We allow $Y$ to be either categorical or continuous, and, for simplicity, assume that $X$ is a random vector in $\real^p$. The proposed work can be readily applied to non-Euclidean $X$, as long as $\H_X$ is specified {\it a priori}. Since the similarity measure of $X$ in (\ref{eq: le original}) is henceforth not involved, we simply denote $w_{Y,ij}$, $\W_Y$, and $\D_Y$ by $w_{ij}$, $\W$, and $\D$, respectively. Following the discussion below (\ref{eq: le}), we assume the symmetry of $w_{ij}$, that is, $w_{ij}=w_{ji}$, without loss of generality. Other notation includes $\Omega_Y$ for the support of $Y$, $\Gs_X$ for the function that maps each $x \in \real^p$ to $(K_X(X_i,x) -  n\inv\sum_{j=1}^n K_X (X_j, x))_{i=1,\ldots,n} \in \real^{n}$, $\Gs_Y$ defined analogously given an RKHS of $Y$ (see Section \ref{sec: sdr} later), $\vecone_n$ for $(1,\ldots, 1)\trans \in \real^n$, and $\zero$ for the origin of a real space where the dimension is omitted if there is no ambiguity.


\section{A new supervised Laplacian embedding method}\label{sec: sle}

Let $U \in \real^d$ be the reduced predictor that generates $U_1,\ldots, U_n$ mentioned in the Introduction. As in \cite{tai2022kernelized}, we assume $U = \phi_0(X)$ where the components of $\phi_0(\cdot)$ belong to a prespecified $\H_X$. Following common practice, we regularise $U$ as follows. First, since translating $U$ by a non-stochastic vector does not affect its predictability for $Y$, we assume without loss of generality that $E\{U\} = E\{\phi_0(X)\} = \zero$. This allows us to restrict $\H_X$ so that any $\phi(\cdot) \in \H_X$ must satisfy $E\{\phi(X)\} = 0$. At the sample level, the restriction modifies the kernel matrix $\K_X$ defined above (\ref{eq: le sup kernel tai22}) to 
\begin{equation}\label{def: GX}
\G_X = \Q \K_X \Q,
\end{equation}
with $\Q$ denoting $\I_n - n\inv\vecone_n \vecone_n\trans$. Second, we normalise $U$ by setting $\var (U) = \I_d$. Recall from the Introduction that $\U$ is defined as $(U_1, \ldots, U_n)\trans$ and the kernelisation implies $\U = \G_X \V$ for some $\V \in \real^{n\times d}$. At the sample level, this normalisation means 
\begin{equation}\label{reg: var U one}
\U\trans \U =  n \I_d \,\, \mbox{ or equivalently } \,\, \V\trans \G_X^2 \V = n \I_d,
\end{equation}
which complies with the convention in unsupervised Laplacian embedding, as seen in (\ref{eq: le}).

Recall that the ineffectiveness of the method in \cite{tai2022kernelized} stems from the inappropriate normalisation constraint (\ref{eq: le sup kernel tai22 constr}), which involves the response $Y$ through $\D$. Naturally, this issue will be addressed if we can modify the similarity measure $w_{ij}$ so that the corresponding $\D$ has constant diagonal entries invariant to $Y$. To achieve this, for any prespecified constant $\delta>0$, we modify each $w_{ij}$ to 
\begin{equation}\label{eq: centred wij}
w_{ij} - \bar w_{i\cdot} - \bar w_{\cdot j} + \bar w_{\cdot\cdot} + \delta, 
\end{equation}
where $\bar w_{i\cdot}$ denotes $n\inv \sum_{j=1}^n w_{ij}$, $\bar w_{\cdot j}$ denotes $n\inv \sum_{i=1}^n w_{ij}$, and $\bar w_{\cdot\cdot}$ denotes $n\inv \sum_{i=1}^n \bar w_{i\cdot}$, that is, the row average, the column average, and the overall average of $\W$, respectively. By simple algebra, this modification transforms $\D$ to $n \delta \I_n$ and consequently turns the normalisation constraint (\ref{eq: le sup kernel tai22 constr}) into 
\begin{equation}\label{eq: normalisation unused}
  \delta \V\trans \G_X^2 \V = \I_d.  
\end{equation}
Since the effectiveness of dimensionality reduction is invariant under any scaling of the reduced predictor, we can always replace $\I_d$ on the right-hand side of (\ref{eq: normalisation unused}) with $\delta \I_d$. This makes (\ref{eq: normalisation unused}) identical to (\ref{reg: var U one}) and therefore unneeded for implementation. In addition, since $\D - \W$ in the objective function (\ref{eq: le sup}) is invariant to $\delta$, we can simply set $\delta = 0$ without changing the results of dimensionality reduction. Denote the modified $w_{ij}$ in (\ref{eq: centred wij}) by $w_{ij}^c$, where the superscript $c$ stands for ``centred'', as $w_{ij}^c$ satisfies $\sum_{i=1}^n w_{ij}^c = \sum_{j=1}^n w_{ij}^c = 0$. 

From a practical point of view, the modification of $w_{ij}$ in (\ref{eq: centred wij}) represents a concern for fairness; that is, we should not treat the importance of the similarity between $Y_i$ and $Y_j$ differently based on their locations in $\Omega_Y$. Specifically, since (\ref{eq: centred wij}) can be written in reverse as $w_{ij} = w_{ij}^c + \bar w_{i\cdot} + \bar w_{\cdot j} - (\bar w_{\cdot\cdot} + \delta)$, under the normalisation constraint (\ref{reg: var U one}), the objective function in (\ref{eq: le sup}) is identical to
\begin{equation}\label{eq: sle uneq wi}
     \sum_{i=1}^n \sum_{j=1}^ n \|U_i - U_j\|^2 w_{ij}^c + 2 \sum_{i=1}^n \sum_{j=1}^ n \|U_i - U_j\|^2 \bar w_{i\cdot}
\end{equation}
up to an additive constant. If $Y_i$ induces a larger $\bar w_{i\cdot}$, then the second term of (\ref{eq: sle uneq wi}) encourages $U_j$ to be closer to $U_i$. Intuitively, for continuous $Y$, this occurs when $Y_i$ is near the centre of the data cloud. Therefore, the weighting mechanism in (\ref{eq: sle uneq wi}) favours better prediction of $U$ for non-outliers in $Y$ than for outliers, making it often incompatible with research interest in practice. A similar discussion for the case of categorical $Y$ can be found in Example \ref{example: delta} below. Setting $\bar w_{i\cdot}$ to be invariant to $Y_i$ represents the situation in which we do not discriminate in the predictability of $U$ across different regions or classes of $Y$, and instead would like $U$ to perform well over the entire population.

\begin{example}\label{example: delta}
Suppose that $Y$ is categorical with $\Omega_Y = \{1,\ldots, K\}$, and its probability mass function is again denoted by $\pi(\cdot)$. We abuse notation by writing $w_{ij}$ as $w_{k,s}$ when $Y_i = k$ and $Y_j = s$, and $\bar w_{i\cdot}$ as $\bar w_{k, \cdot}$ when $Y_i = k$. For simplicity, we set $w_{k,s} = \ind(k=s)$. For all large $n$, $\pi(k) > \pi(s)$ is equivalent to $\bar w_{k, \cdot}> \bar w_{s,\cdot}$; that is, the reduced predictor $U$ that minimises (\ref{eq: sle uneq wi}) must have smaller variation in a larger class of $Y$, which is undesirable unless endorsed by domain knowledge. 
\end{example}

By replacing $w_{ij}$ with $w_{ij}^c$, $\W$ is modified to $\W^c = \Q \W \Q$ and $\D$ is modified to $\zero \in \real^{n\times n}$. Thus, the minimisation of (\ref{eq: le sup kernel tai22}) becomes the minimisation of $\tr\{\V\trans \G_X (-\W^c)\G_X \V\}$, or equivalently the maximisation of $\tr(\V\trans \G_X \W^c \G_X \V)$. Using the normalisation constraint (\ref{reg: var U one}), we propose 
\begin{equation}\label{def: sle 1}
\max_{\V \in \real^{n\times d}}  \tr(\V\trans \G_X \W^c \G_X \V) \quad \mbox{subject to } \V\trans \G_X^2 \V = n \I_d.
\end{equation}
Let $\widehat \V$ be the maximiser in (\ref{def: sle 1}). Our reduced predictor, denoted by $\widehat U$, is $\Gs_X\trans \widehat \V$ where $\Gs_X$ is defined at the end of the Introduction. 
To solve the generalised eigen-decomposition problem (\ref{def: sle 1}), we follow common practice \citep{lee13} to first take $\Theta = n^{-1/2}\G_X \V$, which satisfies $\Theta\trans \Theta= \I_d$, then derive $\widehat \Theta$ as the set of $d$ leading eigenvectors of 
\begin{equation}\label{eq: sle 1 eigen}
\G_X^\dagger \G_X \W^c \G_X \G_X^\dagger, 
\end{equation}
where $\G_X^\dagger$ denotes the Moore-Penrose inverse of $\G_X$, and finally take $\widehat \V$ as $n^{1/2} \G_X^\dagger \widehat \Theta$. To approximate $\G_X^\dagger$ in (\ref{eq: sle 1 eigen}), we use Tikhonov regularisation \citep{groetsch1984, lee13}, that is, $(\G_X + \epsilon_n \I_n)\inv$ where $\epsilon_n$ is a negligible positive scalar such as $(n \log n)\inv$.

In the above implementation procedure, we also need to determine $d$, the dimension of the reduced predictor. In the statistical literature, this is called order determination and has been thoroughly studied in the scenario of linear dimensionality reduction \citep{bai2002, bura2011, luo2021}, where the dimension of the reduced predictor coincides with the number of nonzero eigenvalues of some matrix-valued parameter. For nonlinear dimensionality reduction, order determination is more complex due to the redundancy issue. That is, although the components of the reduced predictor are linearly independent by virtue of the eigen-decomposition of (\ref{eq: sle 1 eigen}), one component may still be a nonlinear function of the others and therefore redundant for modelling $Y$ even if the corresponding eigenvalue of (\ref{eq: sle 1 eigen}) is nonzero. 
This issue can be addressed by imposing probabilistic independence among the components of the reduced predictor \citep{bengio13, huang24}, but at the expense of complicating the implementation and jeopardising interpretability. To avoid distractions, we leave modifications of the proposed method in this direction to future work. In what follows, we provide an order-determination method based on a variable importance score that is more relevant than the eigenvalues of (\ref{eq: sle 1 eigen}).


\newcommand\comp[1]{^{(#1)}}
\def\M{{\mathrm M}}

Suppose $d_\M$ is a prespecified strict upper bound of $d$, and let $(\widehat U\comp 1,\ldots, \widehat U\comp {d_\M})\trans$ be the maximiser in (\ref{def: sle 1}) when $d$ is misspecified as $d_\M$. Because the objective function in (\ref{def: sle 1}) can be decomposed as the sum of $d$ functions, each involving exactly one column of $\V$, the maximiser in (\ref{def: sle 1}) satisfies a nesting property: for any $k < d_\M$, it is exactly $(\widehat U\comp 1,\ldots, \widehat U\comp {k})\trans$ when $d$ is specified as $k$. In particular, $\widehat U$ defined above consists of $\widehat U\comp 1,\ldots, \widehat U\comp {d}$. Given the nesting property, choosing the appropriate dimension of the reduced predictor boils down to deciding whether each $\widehat U\comp {k+1}$ is informative for modelling $Y$ in the presence of the preceding components $(\widehat U\comp 1,\ldots, \widehat U\comp {k})\trans$. To this end, we calculate the sample conditional distance correlation \citep{wang2015} between $\widehat U\comp {k+1}$ and $Y$ given $(\widehat U\comp 1,\ldots, \widehat U\comp {k})$, denoted by $\widehat \tau_k$. When $Y$ is categorical with $K$ categories, we assign these categories numerical values $1,\ldots, K$ to compute $\widehat \tau_k$. Each $\widehat \tau_k$ ranges from $0$ to $1$ and is negligible if and only if $\widehat U \comp {k+1}$ is nearly uninformative for $Y$ in the presence of $(\widehat U\comp 1,\ldots, \widehat U\comp {k})$. To accommodate the extreme case where $d=0$, which means that $X$ is completely uninformative for $Y$, we additionally set $\widehat \tau_0$ to the sample marginal distance correlation between $\widehat U\comp 1$ and $Y$. Following the spirit of the Bayesian information criterion \citep[BIC;][]{schwarz1978}, we determine $d$ as the minimiser of
\begin{equation}\label{eq: BIC}
\max\{\widehat \tau_\ell: \ell=k,\ldots,d_\M-1\} + \rho_n k
\end{equation}
over $k \in \{0,\ldots, d_\M - 1\}$, where $\rho_n$ is a positive scalar that slowly vanishes as $n\rightarrow \infty$. In Section \ref{sec: sle theory}, we will revisit this BIC-type criterion and justify its consistency.

We conclude this section by drawing attention to a slight variation of the proposed method. Since $\W^c$ yields $\D=\zero$, the normalisation constraint (\ref{reg: var U one}) is not needed for the equality of the objective functions in (\ref{eq: le sup kernel tai22}) and (\ref{def: sle 1}). Thus, we may use $\V\trans \G_X \V = n \I_d$, another commonly used normalisation in the literature on RKHS-based methods, in place of $\V\trans \G_X^2 \V = n \I_d$ in (\ref{def: sle 1}). Let $\widetilde \V$ be the maximiser in (\ref{def: sle 1}) under this normalisation. Compared to $\widehat \V$ above, the implementation of $\widetilde \V$ is slightly easier, as it has a closed form without involving $\G_X^\dagger$: let $\widetilde \Theta$ be the set of $d$ leading eigenvectors of $\G_X^{1/2} \W^c \G_X^{1/2}$ and $\widetilde \Lambda$ be the diagonal matrix whose diagonal entries are the corresponding eigenvalues, so that $\G_X^{1/2} \W^c \G_X^{1/2} \widetilde \Theta = \widetilde \Theta \widetilde \Lambda$; then $\widetilde \V = n^{1/2} \W^c \G_X^{1/2} \widetilde \Theta \widetilde \Lambda\inv$. Our numerical experience suggests that $\Gs_X\trans \widetilde \V$ performs similarly to $\Gs_X\trans \widehat \V$, so we omit the details of this variant in the rest of the article.

\section{Effectiveness of supervised Laplacian embedding}\label{sec: sle theory}

We now justify the proposed supervised Laplacian embedding method at the population level, that is, when it delivers an effective reduced predictor preserving the information in $X$ for modelling $Y$. Recall from the beginning of Section \ref{sec: sle} that the reduced predictor at the population level is denoted by $\phi_0(X)$, and its sample estimator is $\Gs_X\trans \widehat \V$ derived from (\ref{def: sle 1}), where the objective function originates from (\ref{eq: le sup}). To ease the presentation, we use $(X, Y)$ and $(X^*, Y^*)$ instead of $(X_i, Y_i)$ and $(X_j,Y_j)$ to denote two random observations. At the population level, a similarity measure of $Y$ refers to a symmetric function $w: \Omega_Y \times \Omega_Y \rightarrow \real$ with $w(y, y^*) = w(y^*, y)$. Its centred version is
\begin{equation*}
    w^c(Y, Y^*) = w(Y, Y^*) - E\{w(Y,Y^*)\mid Y\} - E\{w(Y,Y^*)\mid Y^*\} + E\{w(Y,Y^*)\},
\end{equation*} 
which satisfies $E\{w^c(Y,Y^*)\mid Y\} = E\{w^c(Y,Y^*)\mid Y^*\} = 0$.
For simplicity, we assume that $w(Y,Y^*)$ is always nonnegative and has bounded support, and we also extend the search range for $\phi_0(\cdot)$ from the normalised function sets in $\H_X$ to those in $\Ls_2(X)$, that is, the collection of measurable functions $\phi(\cdot): \real^p \rightarrow \real^d$ satisfying $E\{\phi(X)\} = \zero$ and $\var\{\phi(X)\} = \I_d$, denoted by $\Ldx d$. Together, these imply that the  population-level counterpart of the objective function in (\ref{eq: le sup}) is 
\begin{equation}\label{def: sle 1 pop}
     E\{\|\phi(X) - \phi(X^*)\|^2 w^c(Y, Y^*)\},
\end{equation}
and $\phi_0(\cdot)$ minimises (\ref{def: sle 1 pop}) over $\phi(\cdot) \in \Ldx d$. For ease of understanding, one may tentatively assume the uniqueness of $\phi_0(X)$ up to an orthogonal linear transformation. A detailed discussion is deferred to Condition (C2) later. Rigorously, uniqueness is also up to values on a null set, that is, with zero probability. To simplify the presentation, we henceforth ignore any differences on such sets.

To define a desired reduced predictor, we assume a low-dimensional structure on the data; that is, there exists $f_0: \real^p \rightarrow \real^{d_0}$ with $d_0 < p$ such that $f_0 (X)$ is sufficient for modelling $Y$ in that 
\begin{equation}\label{asm: sdr}
    Y \indep X \mid f_0(X),
\end{equation}
where $\indep$ denotes probabilistic independence. The identifiability of $f_0(X)$ in (\ref{asm: sdr}) is only up to an invertible transformation, which boils down to the $\sigma$-field that $f_0(X)$ generates. Under fairly general conditions, \cite{lee13} proved the unique existence of the smallest $\sigma$-field, called the central $\sigma$-field and denoted by $\sigma_{Y|X}$, whose arbitrary generating set satisfies (\ref{asm: sdr}). Without loss of generality, we assume that $f_0(X)$ is minimal in the sense that it generates $\sigma_{Y|X}$. Under this assumption, $f_0(X)$ serves as the benchmark, up to an invertible transformation, for any dimensionality reduction procedure. We next study when the proposed reduced predictor $\phi_0 (X)$ is a transformation of $f_0 (X)$, and when this transformation is invertible. From (\ref{def: sle 1 pop}), the answer to these questions hinges on an appropriate choice of $w(Y, Y^*)$ or equivalently $w^c(Y, Y^*)$.

As mentioned in the Introduction, although it has been widely used in the literature on Laplacian embedding and general dimensionality reduction, the similarity measure $w(Y,Y^*)$ is an intuitive concept that has not been rigorously formulated. The following condition provides a novel regulation on $w(Y,Y^*)$, and can be regarded as its rigorous definition. Let $\Ls_2^c(Y)$ be the subcollection of $\Ls_2(Y)$ consisting of centred $\eta(Y)$ with $E\{\eta(Y)\} = 0$. 

\vspace{0.2cm}

\noindent
(C1) For any nondegenerate $\eta(Y) \in \Ls_2^c(Y)$, we always have $E\{\eta(Y) \eta(Y^*) w (Y, Y^*)\} > 0$.

\vspace{0.2cm}

\noindent
In terms of $w^c (Y, Y^*)$, Condition (C1) can be restated as: for any nondegenerate $\eta(Y) \in \Ls_2(Y)$,
\begin{equation}\label{asm: C1 equiv}
E\{\eta(Y) \eta(Y^*) \, w^c (Y, Y^*)\} > 0. 
\end{equation}
When $Y$ is categorical, this means that $(w^c(k,s))_{k,s \in \Omega_Y}$ is a positive semi-definite matrix with the null space spanned by $\vecone$. At the sample level, Condition (C1) implies that for any nondegenerate $\eta(Y) \in \Ls_2(Y)$,
$\sum_{i=1}^n \sum_{j=1}^n \eta(Y_i) \eta(Y_j) w^c_{ij} > 0$
for all large $n$; that is, $\W^c$ induced by $w(Y,Y^*)$ is asymptotically positive semi-definite with the null space spanned by $\vecone_n$. 

To interpret Condition (C1), suppose that $Y$ is continuous, and denote the probability density function of $Y$ by $\zeta(\cdot)$. For any nondegenerate $\eta(Y) \in \Ls_2^c (Y)$, the independence between $Y$ and $Y^*$ means that $\eta(Y)$ and $\eta(Y^*)$ are uncorrelated. In contrast, since $w(Y,Y^*)$ is nonnegative and bounded, there exists some constant $a>0$ such that $a\inv \zeta(y) \zeta(y^*) w(y,y^*)$ is a probability density function, under which (C1) implies a positive correlation between $\eta(Y)$ and $\eta(Y^*)$. Thus, (C1) can be understood as that, by changing the joint density of $(Y, Y^*)$ from $\zeta(y) \zeta(y^*)$ to $a\inv \zeta(y) \zeta(y^*) w(y,y^*)$, $Y$ and $Y^*$ become similar to each other in all respects. In this sense, (C1) is also a natural definition of the similarity measure $w (Y, Y^*)$. The same interpretation holds when $Y$ is categorical.

Recall from the Introduction that in the literature on dimensionality reduction, the similarity measure $w(Y,Y^*)$ refers to a function that reflects the metric of $Y$. Thus, it should decrease monotonically with the distance (e.g., Euclidean) between $Y$ and $Y^*$, which requires that
\begin{equation}\label{def: sm literature}
w(y,y) > w(y,y^*) \mbox{ for any } y, y^* \in \Omega_Y \mbox{ with } y\neq y^*. 
\end{equation}
Condition (C1) does not contradict this intuition for commonly used similarity measures in practice, such as $\exp \{-(Y - Y^*)^2/(2 \sigma^2)\}$ where $\sigma^2$ is a positive constant. As will be seen in Section \ref{sec: sdr} later, the reason is that these similarity measures are in fact kernel functions of certain RKHSs, which automatically satisfy (C1). In particular, a family of similarity measures in the form of $w(Y - Y^*)$ for continuous $Y$ is closely related to a concept known as a continuous negative definite function in the statistical literature \citep{li2023}. Nonetheless, (C1) and (\ref{def: sm literature}) still differ conceptually. Interestingly, the difference is two-sided: while (C1) excludes some similarity measures that satisfy (\ref{def: sm literature}), it also includes others that violate (\ref{def: sm literature}). We next explore the details.

\begin{figure}[h]
     \centering
            \includegraphics[page=1, width=0.36\textwidth]{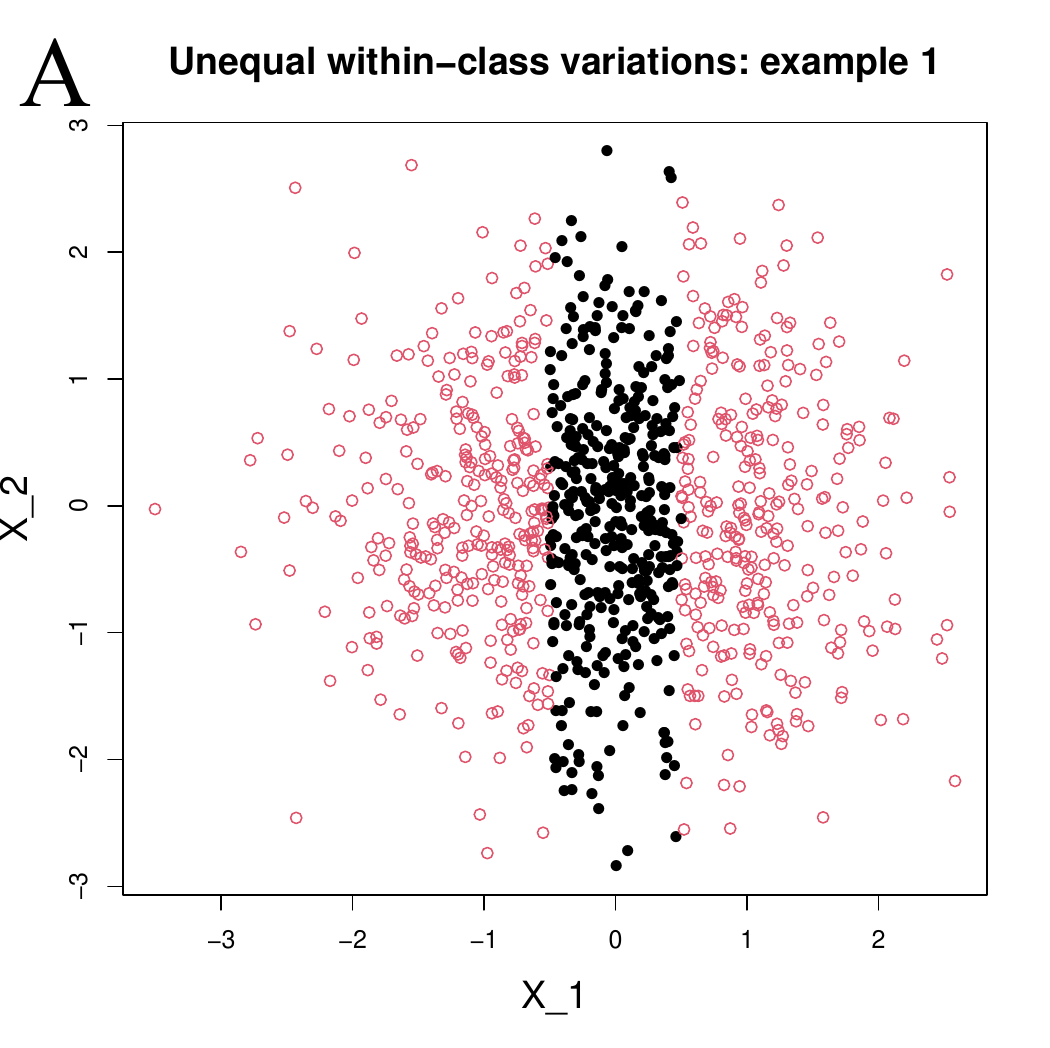}
            \includegraphics[page=2, width=0.36\textwidth]{uneq_wii.pdf}
    \caption{The scatter plots of $(X\comp 1, X\comp 2)$, with ``$\cdot$'' denoting the class $Y=1$ and ``$\circ$'' denoting the class $Y=2$. In Panel A, $X$ follows $N(\zero, \I_2)$ and $Y=\ind (|X\comp 1| > 0.5) + 1$; in Panel B, $X$ is set as $(Y-1+\varepsilon_R) \varepsilon_U$, where $\varepsilon_R$ is uniformly distributed on $(0, 0.1)$, $\varepsilon_U$ is uniformly distributed on the unit sphere in $\real^2$, and both $\varepsilon_R$ and $\varepsilon_U$ are independent of $Y$.}
    \label{fig: uneq wii}
\end{figure}

On the one hand, Condition (C1) allows $\min\{w(y,y), w(y^*,y^*)\} < w (y,y^*)$ when $y$ differs from $y^*$, thereby relaxing (\ref{def: sm literature}). In particular, if $Y$ is binary with $\Omega_Y=\{1, 2\}$, then (C1) holds as long as $w_{1,1} w_{2,2} > w_{1,2}^2$, which is strictly weaker than $\min\{w_{1,1}, w_{2,2}\} > w_{1,2}$ required by (\ref{def: sm literature}). This flexibility is crucial in at least two cases. First, if $w_{1,1} > w_{2,2}$, which is sensible when the true classification of $Y=1$ (such as a true positive) is more important than that of $Y=2$ (such as a true negative), then allowing $w_{1,2}>w_{2,2}$ will enhance the freedom in tuning $w_{1,2}$ and thereby improve the effectiveness of dimensionality reduction. Second, as illustrated in both panels of Figure \ref{fig: uneq wii}, if subjects from different classes of $Y$ are, on average, more similar to each other than those from the class of $Y=2$, then it is naturally desirable to set $w_{1,2}>w_{2,2}$.  

On the other hand, a similarity measure that satisfies (\ref{def: sm literature}) may violate Condition (C1), even if it decreases strictly with the distance between $y$ and $y^*$. When $Y$ is continuous, this follows from results on conditional negative definite functions \citep[Section 3.3,][]{berg1984harmonic}, which imply that $-\|y-y^*\|^c$ violates (C1) when $c>2$. As another example, suppose $Y$ is categorical with $\Omega_Y = \{1, 2, 3\}$ and 
\begin{equation}\label{eq: sm case 1}
(w_{k,s})_{k,s=1,2,3} = \begin{pmatrix} 1 & \delta & 0 \\ \delta & 1 & \delta \\ 0 & \delta & 1 \end{pmatrix},
\end{equation}
then $w(\cdot,\cdot)$ satisfies (\ref{def: sm literature}) as long as $\delta \in (0, 1)$. However, if $\delta >1/\sqrt 2$, then $w(\cdot,\cdot)$ violates Condition (C1); that is, the corresponding $(w_{k,s}^c)_{k,s=1,2,3}$ is no longer positive semi-definite. 


Note that if $\delta >1/\sqrt 2$ in (\ref{eq: sm case 1}), the proposed supervised Laplacian embedding method may not deliver the desired reduced predictor. To see this, suppose $X = (X\comp 1,X\comp 2, X\comp 3)\trans$ follows $N(\zero, \I_3)$,
\begin{equation}\label{eq: counter example}
Y = \left\{ \begin{matrix} 1 & \mbox{ if } X\comp 1 < 0 \mbox{ and } X \comp 2 < 0, \\
2 & \mbox{ if } X\comp 1 > 0, \\
3 & \mbox{ if } X\comp 1 < 0 \mbox{ and } X\comp 2 > 0.
\end{matrix}\right. 
\end{equation}
and we restrict the search range for the reduced predictor to linear transformations of $X$. Then (\ref{def: sle 1 pop}) is minimised by $\pm X\comp 2$ when $d=1$, and by $(X\comp 2,X\comp 3)$ up to an orthogonal linear transformation when $d=2$, neither of which captures $X\comp 1$ and distinguishes the second class from the other two classes. The same limitation applies to the existing supervised Laplacian embedding methods mentioned in the Introduction. In light of (\ref{eq: le sup}), the failure of supervised Laplacian embedding in this case stems from the nature of (\ref{eq: sm case 1}), which excessively emphasises the difference between the first and third classes over those between the other classes: by (\ref{eq: sm case 1}) and (\ref{eq: counter example}), $X\comp 1_i$ is similar to $X\comp 1_j$ when $(Y_i,Y_j) = (1, 3)$ or $(3, 1)$, in which case $w_{ij}$ is small; $X\comp 1_i$ is dissimilar to $X\comp 1_j$ when $(Y_i, Y_j) = (1,2)$, $(2,1)$, $(2,3)$, or $(3,2)$, in which case $w_{ij}$ is large. Since a minimiser of (\ref{def: sle 1 pop}) must exhibit the opposite pattern, it is irrelevant to $X\comp 1$. 



From the discussion above, Condition (C1) is necessary for the effectiveness of supervised Laplacian embedding. We next justify its sufficiency. First, we formalise the uniqueness of $\phi_0 (X)$ mentioned below (\ref{def: sle 1 pop}) by regulating the behaviour of (\ref{def: sle 1 pop}) on $\Ldfx d$, the subcollection of $\Ldx d$ consisting of centred and normalised transformations of $f_0(X)$. 

\vspace{0.2cm}

\noindent
(C2) For any $d \geq 1$,  the minimiser of (\ref{def: sle 1 pop}) over $\phi(X) = g\circ f_0(X) \in \Ldfx d$ is always unique up to an orthogonal linear transformation. 

\vspace{0.2cm}

\noindent
Because $f_0 (X)$ represents the informative part of $X$ for modelling $Y$, it reduces to the informative linear transformations of $X$ in the scenario of linear dimensionality reduction. Condition (C2) therefore resembles the assumption often adopted therein \citep{fan2017}, which requires certain symmetric matrix-valued parameters to have distinct nonzero eigenvalues. To accommodate the greater flexibility of $\phi(X)$ caused by its nonlinearity, here we no longer restrict $d$ to be a small integer as in the scenario of linear dimensionality reduction. 
The following theorem shows that, under Conditions (C1) and (C2), $\phi_0(X)$ is a well-defined transformation of $f_0 (X)$, that is, it partially recovers the desired reduced predictor without introducing extraneous noise for subsequent analysis.

\begin{theorem}\label{thm: effectiveness}
Under Conditions (C1) and (C2), for any $d>0$, the minimiser $\phi_0(X)$ of (\ref{def: sle 1 pop}) over $\Ldx d$ is unique up to an orthogonal linear transformation, and must be a transformation of $f_0 (X)$. 
\end{theorem}

To ensure that $\phi_0(X)$ fully recovers the desired reduced predictor $f_0 (X)$, which is a stronger result than Theorem \ref{thm: effectiveness}, we need to additionally assume that any effect of $f_0(X)$ on $Y$ can be captured by the objective function (\ref{def: sle 1 pop}). From the proof of Theorem \ref{thm: effectiveness} (see the Supplementary Material), the link connecting the minimiser of (\ref{def: sle 1 pop}), Condition (C1), and the desired $f_0 (X)$ defined in (\ref{asm: sdr}) is 
\begin{equation}\label{eq: p.s.d.}
E\{\phi(X) \mid Y\} = E[E\{\phi(X)\mid f_0 (X)\}\mid Y] \mbox{ for any } \phi(X) \in \Ls_2(X). 
\end{equation}
If there exists $g \circ f_0 (X) \in \Ldfx 1$ such that $E\{g\circ f_0 (X) \mid Y\} = 0$, then $g \circ f_0 (X)$ will yield the value zero in (\ref{eq: p.s.d.}) as if it were independent of $Y$, and is therefore undetectable by minimising (\ref{def: sle 1 pop}). Thus, for $\phi_0(X)$ to fully recover $f_0(X)$, it is necessary to assume that the collection of $g\circ f_0 (X)$ with $E\{g\circ f_0 (X) \mid Y\} = 0$ is not so large as to mask a sub-$\sigma$-field of $f_0(X)$, that is,

\vspace{0.2cm}

\noindent
(C3) The subcollection of $g\circ f_0(X) \in \Ldfx 1$ with nondegenerate $E\{g\circ f_0 (X) \mid Y\}$ forms an invertible transformation of $f_0(X)$.

\vspace{0.2cm}

\noindent
A strengthened condition is to preclude any $g \circ f_0 (X) \in \Ldfx 1$ with degenerate $E\{g\circ f_0 (X) \mid Y\}$, that is, for any $g \circ f_0(X) \in \Ls_2 (f_0(X))$, $E\{g\circ f_0 (X) \mid Y\} = 0$ implies $g \circ f_0 (X) = 0$. 
This was adopted in \cite{lee13}, where it was called the completeness assumption and its generality was justified in the scenarios of both homoscedastic regression (see their Proposition 1) and inverse regression (see their Proposition 2). In this regard, Condition (C3) is fairly general. 


Under Condition (C3), we can expect $\phi_0(X)$ to fully recover $f_0(X)$ for an appropriate choice of $d$. Because we need $\phi_0(X)$ to be lower-dimensional than $X$ itself, we additionally require that the dominant signals in $\Ldfx 1$ detected by (\ref{def: sle 1 pop}) form an invertible transformation of $f_0(X)$. This is formalised in the following regularity condition, which strengthens both (C2) and (C3). 
Due to the possible redundancy of $\phi_0(X)$ mentioned below (\ref{eq: sle 1 eigen}), this condition is also more restrictive than the general singularity assumption (\ref{asm: sdr})---a price paid for algorithmic simplicity via eigen-decomposition.

\vspace{0.2cm}

\noindent
(C4) For some $d < p$,  the unique minimiser of (\ref{def: sle 1 pop}) over $\phi(X) = g\circ f_0(X) \in \Ldfx d$ is an invertible transformation of $f_0 (X)$. 

\vspace{0.2cm}

\noindent
Hereafter, we take $d$ to be the minimal integer that satisfies Condition (C4) and regard it as the dimension of the desired $\phi_0(X)$. Under Condition (C4), the full recovery of $f_0(X)$ by $\phi_0(X)$ is justified in the following theorem, which summarises the discussion above. 

\begin{theorem}\label{thm: exhaustiveness}
Under Conditions (C1)--(C4), the minimiser $\phi_0(X)$ of (\ref{def: sle 1 pop}) over $\Ldx d$ must be an invertible transformation of $f_0 (X)$. 
\end{theorem}

By Theorem \ref{thm: exhaustiveness}, the proposed supervised Laplacian embedding method provides an effective reduced predictor for subsequent modelling. In practice, this also relies on the negligibility of the sampling error and a consistent selection of $d$. Under certain regularity conditions on $\H_X$, the former is guaranteed with the sampling error being $\bop (n^{-1/2})$; see the relevant asymptotic studies of RKHS-based methods \citep{li2017,li2024}. To select $d$, recall that we proposed a BIC-type criterion in Section \ref{sec: sle}. We next briefly justify the selection consistency of this criterion. 

Let $(U\comp 1,\ldots, U\comp {d_\M})\trans$ be the minimiser of (\ref{def: sle 1 pop}) over $\Ldfx {d_\M}$. As with (\ref{def: sle 1}), the objective function in (\ref{def: sle 1 pop}) is a sum of $d$ functions, each of which involves one component of $\phi(X)$. Thus, $(U\comp 1,\ldots, U\comp {d_\M})\trans$ satisfies the nesting property that, for any $k< d_\M$, its first $k$ components minimises (\ref{def: sle 1 pop}) when $d$ is set to $k$. In particular, $\phi_0(X)$ is identical to $(U\comp 1,\ldots, U\comp {d})\trans$. By Theorem \ref{thm: exhaustiveness}, these imply that  $d$ is the minimal value of $k$ for which $X\indep Y \mid (U\comp 1,\ldots, U\comp {k})\trans$ or equivalently 
\begin{equation}\label{eq: uk indep y u1tok-1}
    U\comp {\ell+1} \indep Y \mid (U\comp 1,\ldots, U\comp {\ell})
\end{equation}
for all $\ell \in \{k,\ldots, d_\M-1\}$. Let $\tau_\ell$ be the conditional distance correlation \citep{wang2015} between $U\comp {\ell+1}$ and $Y$ given $(U\comp 1,\ldots, U\comp {\ell})$, which is the population limit of $\widehat \tau_\ell$ used in (\ref{eq: BIC}) for the proposed BIC-type criterion. Since $\tau_\ell$, which is always nonnegative, reaches zero if and only if (\ref{eq: uk indep y u1tok-1}) holds, $\max\{\tau_\ell: \ell=k,\ldots,d_\M-1\}$ conveys a scree plot that is positive when $k < d$ and reaches zero when $k \geq d$. Consequently, the minimiser of (\ref{eq: BIC}) converges to $d$ as long as $\rho_n$ vanishes at a slower rate than $\max\{|\widehat \tau_\ell|: \ell = d,\ldots, d_\M-1\}$, which is $\bop (n^{-1/2})$ under appropriate regularity conditions similarly to the sampling error mentioned above. The proof resembles that of the classical BIC and is therefore omitted. In practice, we set $d_\M = \min\{10, p\}$ and $\rho_n = 0.5 \, n^{-1/3}$.


From the discussion above, Condition (C1) on the choice of the similarity measure is the foundation of the effectiveness of the proposed supervised Laplacian embedding method. This condition also connects supervised Laplacian embedding to the rich literature on sufficient dimension reduction in statistical research, as we explore in the next section.



\section{Connection with sufficient dimension reduction} \label{sec: sdr}

\def\lspf{\Ls_2^c \{X \ominus  f_0(X)\}}
\def\cran{\overline{\mathrm{ran}}}

In the statistical literature, dimensionality reduction is commonly called dimension reduction. In supervised learning, it is referred to as sufficient dimension reduction \citep[SDR;][]{li1989, li1991} or, more specifically, nonlinear SDR \citep{lee13,sdr18} when a nonlinear transformation of $X$ is allowed to generate the reduced predictor. The term ``sufficiency'' refers to (\ref{asm: sdr}) in Section \ref{sec: sle theory}; that is, the reduced predictor $f_0 (X)$ preserves the information in $X$ for modelling $Y$.

Unlike dimensionality reduction where mainstream research focuses on the alignment between $f_0 (X)$ and $Y$, existing nonlinear SDR methods commonly recover $f_0 (X)$ by characterising the orthogonal complement of $\Ls_2^c\{f_0(X)\}$ in $\Ls_2^c(X)$, denoted by $\lspf$, which is the collection of $\alpha(X) \in \Ls_2^c(X)$ such that $E\{\alpha(X) \mid f_0 (X)\} = 0$. Here, $\Ls_2^c(X)$ and $\Ls_2^c\{f_0(X)\}$ respectively denote the subcollections of $\Ls_2(X)$ and $\Ls_2\{f_0(X)\}$ consisting of centred elements with mean zero, analogously to $\Ls_2^c (Y)$ defined above Condition (C1). For example, a state-of-the-art nonlinear SDR method, called generalised sliced inverse regression \citep[GSIR;][]{lee13}, uses (\ref{eq: p.s.d.}) to characterise each $\alpha(X) \in \lspf$ by
\begin{equation}\label{eq: gsir pop}
E\{\alpha(X)\mid Y \} = 0. 
\end{equation}
That is, $\alpha(X)$ must fall in the kernel of the conditional mean operator $E_{X|Y}: \Ls_2^c (X) \rightarrow \Ls_2^c (Y)$, which maps each $\phi(X)$ to $E\{\phi(X) \mid Y\}$. This induces a collection of transformations of $f_0 (X)$ identified with the range of $E^*_{X|Y} E_{X|Y}: \Ls_2^c (X) \rightarrow \Ls_2^c (X)$, where $E_{X|Y}^*$ denotes the adjoint operator of $E_{X|Y}$. For computational convenience, \cite{lee13} approximated $\Ls_2^c(X)$ by a suitable RKHS of $X$ and approximated $\Ls_2^c (Y)$ by a suitable RKHS of $Y$ (if $Y$ is continuous), denoted by $\H_X$ as in Section \ref{sec: sle} above and by $\H_Y$, respectively. Accordingly, $E_{X|Y}$ is modified to a (linear) operator from $\H_X$ to $\H_Y$, and the range of $E^*_{X|Y} E_{X|Y}: \H_X \rightarrow \H_X$ is modified to its closure under the $L_2$ norm. More details about these modifications can be found in \cite{lee13} and \cite{li1991}.



In multiple respects, the proposed supervised Laplacian embedding method generalises GSIR. To see this, we reformulate GSIR in a similar manner to (\ref{def: sle 1 pop}). Let $K_Y(\cdot,\cdot)$ be the kernel function of $\H_Y$ and let $\{\psi_i(\cdot): i = 1,\ldots, q\}$ be its eigenfunctions associated with positive eigenvalues, where $q$ can be infinite and each $\psi_i(\cdot)$ satisfies $E\{\psi_i(Y)\} = 0$. 

\begin{theorem}\label{thm: gsir reformulation}
The range of $E^*_{X|Y} E_{X|Y}: \H_X \rightarrow \H_X$, regarded as a linear space of functions, is spanned by the eigenfunctions of the negative semi-definite mapping
\begin{equation}\label{eq: gsir pop equiv}
(\phi(\cdot) \in \H_X, \widetilde \phi(\cdot)\in \H_X) \mapsto E[\{\phi(X) - \widetilde \phi(X^*)\}^2 \wgsir(Y,Y^*)]
\end{equation}
that correspond to negative eigenvalues, where $\wgsir(Y,Y^*)$ denotes $\sum_{i=1}^ q \psi_i(Y) \psi_i(Y^*)$ and $(X^*,Y^*)$ again denotes an independent copy of $(X, Y)$.    
\end{theorem}

By Theorem \ref{thm: gsir reformulation}, GSIR equivalently seeks the unique minimiser of 
\begin{equation}\label{eq: gsir pop eigen}
E\{\|\Phi(X) - \Phi(X^*)\|^2 \wgsir(Y,Y^*)\}
\end{equation}
over all $\Phi(\cdot)\subset \H_X$ consisting of $r$ elements with $\var\{\Phi(X)\} = \I_r$, where 
$r$ denotes the rank of (\ref{eq: gsir pop equiv}) and is assumed to be less than $p$ for the effectiveness of GSIR. By simple algebra, $\wgsir(Y,Y^*)$ is centred in the sense that $E\{\wgsir(Y,Y^*) \mid Y\} = E\{\wgsir(Y,Y^*) \mid Y^*\} = 0$, and satisfies Condition (C1) in the form of (\ref{asm: C1 equiv}). Therefore, it is a valid choice of the centred similarity measure $w^c (Y,Y^*)$. In particular, when $Y$ is categorical, $\wgsir(Y,Y^*)$ coincides with the centred version of $\ind(Y=Y^*)$, a similarity measure commonly used in classification. 

Comparing (\ref{eq: gsir pop eigen}) with (\ref{def: sle 1 pop}), GSIR is essentially a special case of the proposed supervised Laplacian embedding method, differing only in that it collects all transformations of $X$ that are orthogonal to those yielding the value zero in (\ref{eq: gsir pop eigen}), whereas we retain only the ``best'' of these transformations. An advantage of the latter is that we allow $r$, the rank of (\ref{eq: gsir pop equiv}), to be infinite, in contrast to the restrictive assumption $r < p$ required for GSIR. In line with these, the consistency of the proposed method discussed in Section \ref{sec: sle theory} parallels that of GSIR, and Condition (C3) again relaxes the completeness assumption required for GSIR to recover an invertible transformation of $f_0(X)$. At the sample level, with $\G_X$ and $\G_Y$ being the kernel matrices corresponding to (centred) $\H_X$ and $\H_Y$, respectively, $E_{X|Y}$ maps each $\Gs_X\trans v$ for $v \in \real^n$ to $\Gs_Y\trans \G_Y^\dagger \G_X v$. Thus, GSIR delivers the reduced predictor $\Gs_X\trans \V$ where $\V \in \real^{n\times r}$ maximises 
\begin{equation}\label{eq: gsir sample}
\tr(\V\trans \G_X \G_Y^\dagger \G_Y^2 \G_Y^\dagger \G_X \V) \quad \mbox{subject to } \V\trans \G_X^2 \V = \I_r.
\end{equation} 
Here, $\G_Y^\dagger$ is approximated by $(\G_Y + \epsilon_n^*\I_n)\inv$ for some prespecified $\epsilon_n^* > 0$, as with the approximation of $\G_X^\dagger$ mentioned below (\ref{eq: sle 1 eigen}). The connection between (\ref{eq: gsir sample}) and (\ref{def: sle 1}) mirrors that between the two methods at the population level described above.





Following a similar idea to GSIR, \cite{chen2024} proposed another family of nonlinear SDR methods that characterise (\ref{eq: gsir pop}) through generalised martingale difference divergence \citep[GMDD;][]{li2023}. That is, let $\Theta(Y - Y^*)$ be a continuous negative definite function taking values in $\real$ such that, for any $y_1, \ldots, y_n \in \Omega_Y$, the $n\times n$ symmetric matrix $(\Theta(y_i - y_j))_{i,j=1, \ldots, n}$ is negative semi-definite with the null space spanned by $\vecone_n$ \citep{li2023}. By construction, $E[\{\alpha(X) - \alpha (X^*)\}^2 \{-\Theta (Y - Y^*)\}]$ is maximised over $\Ldx 1$ when $\alpha (X)$ satisfies (\ref{eq: gsir pop}), which holds if $\alpha (X)$ belongs to $\lspf$. Since all transformations of $f_0 (X)$ are orthogonal to $\lspf$, there exists an appropriate integer $r>0$ such that the minimiser of 
\begin{equation}\label{eq: gmdd}
   E\left[\|\Phi(X) - \Phi(X^*)\|^2 \{-\Theta(Y -Y^*)\}\right] 
\end{equation}
over $\Phi(X) \in \Ldx r$ is a transformation of $f_0 (X)$. This minimiser is the reduced predictor of the GMDD-based nonlinear SDR methods at the population level.

Comparing (\ref{eq: gmdd}) with (\ref{def: sle 1 pop}), the GMDD-based nonlinear SDR methods can also be regarded as a special case of the proposed supervised Laplacian embedding method, in which the centred similarity measure $w^c(Y, Y^*)$ is taken to be $-\Theta(Y -Y^*)$. Unlike $\wgsir(Y,Y^*)$ used in GSIR, $-\Theta(Y -Y^*)$ is based directly on the distance between $Y$ and $Y^*$, and its positive semi-definiteness mentioned above coincides with Condition (C1), again justifying the necessity of the latter. The construction of a valid $\Theta(Y-Y^*)$ has been widely studied in the literature with the aid of the L{\'e}vy measure and the L{\'e}vy--Khintchine formula \citep{van2012, bottcher2013, bottcher2018}, providing a guideline for constructing $w^c(Y,Y^*)$ that satisfies Condition (C1). See some examples along with the associated L{\'e}vy measure in Table 1 of \cite{chen2024}, where all functions decrease monotonically with the distance between $Y$ and $Y^*$. In \cite{chen2024}, $f_0 (X)$ is assumed to lie in a deep neural network class rather than in an RKHS, but this does not alter the nonlinear SDR theory at the population level.

To the best of our knowledge, the discussion above is the first time in the literature that bridges the gap between machine learning research on dimensionality reduction and statistical research on nonlinear SDR. As seen above, the former focuses on preserving the similarity among the observations of $Y$, while the latter focuses on characterising the conditional independence relationship (\ref{asm: sdr}) that reflects the sufficiency of the reduced predictor. The connection uncovered here offers an alternative perspective on both sides and sheds light on their future development by drawing on each other's literature. We outline several directions below, of which the first two concern the development of nonlinear SDR, and the last two concern the development of dimensionality reduction. An example of the latter is elaborated in Section \ref{sec: second-order sle}. A full investigation of these directions requires tremendous effort and is deferred to future work.

First, we now have an algorithmic perspective on GSIR, GMDD-based methods, and potentially other nonlinear SDR methods. That is, instead of focusing on the orthogonal complement of $f_0 (X)$ characterised by (\ref{eq: gsir pop}), we can regard these methods as directly recovering a transformation of $X$ that best captures the similarity among the observations of $Y$. This improves the interpretability of the reduced predictor obtained from nonlinear SDR, and can offer practical guidance for choosing an appropriate $ \H_Y$ based on whether the associated similarity measure suits the nature of the data or the research interest. From the comparison between the proposed supervised Laplacian embedding method and GSIR above, the algorithmic perspective also permits the use of a lower-dimensional reduced predictor than the original nonlinear SDR methods. This enhances the parsimony of the reduced predictor and also relaxes the singularity assumption on the data, thereby facilitating the effectiveness and applicability of nonlinear SDR in practice. 

Second, we can generate more flexible nonlinear SDR methods by, for example, replacing $\wgsir$ in (\ref{eq: gsir pop equiv}) with a more general similarity measure such as the centred version of $K_Y(Y,Y^*)$ given $\H_Y$, which simplifies the computation by avoiding the use of $\G_Y^\dagger$. In light of the discussion below (\ref{eq: gsir pop}), we speculate that this generalisation corresponds to recovering the range of $E_{X|Y}^* R_{Y} E_{X|Y}$, where $R_{Y}$ is an arbitrary operator mapping the range of $E_{X|Y}$ onto itself. The consistency of nonlinear SDR methods under different choices of $R_Y$ also partially explains the commonly observed robustness of GSIR to the approximation of $\G_Y^\dagger$ mentioned below (\ref{eq: gsir sample}); that is, any $\epsilon_n^*>0$ appears to yield a similar reduced predictor regardless of whether the corresponding $(\G_Y + \epsilon_n^* \I_n)\inv$ actually approximates $\G_Y^\dagger$. More discussion is deferred to the numerical studies in Section \ref{sec: simulation} later. To generalise the GMDD-based methods, we may use more complex functions of $(Y,Y^*)$ beyond those depending only on $Y-Y^*$, and a similar theory applies.


Third, we now have a theoretical perspective on supervised Laplacian embedding, which, as seen above, provides suitable regulations that guarantee the sufficiency of the reduced predictor. Following the modification from unsupervised to supervised Laplacian embedding, a variety of other unsupervised dimensionality reduction methods, such as locality preserving projections \citep{he2003LPP}, t-distributed stochastic neighbour embedding \citep[t-SNE;][]{balamurali2023}, and uniform manifold approximation and projection \citep{mcinnes2018umap}, can be adapted into supervised methods centred on aligning the reduced predictor with the response in a certain sense. For example, \cite{cheng2021} proposed a supervised t-SNE that minimises the KL-divergence between the reduced predictor and the response. As with the proposed work, SDR theory may also be applied to regulate these methods and justify their effectiveness. In particular, we speculate that the supervised t-SNE proposed in \cite{cheng2021} can be regulated using the linear SDR theory established in \cite{yin2005direction} and \cite{yin2008successive}, which is also based on minimising KL-divergence. 

Finally, in the literature on nonlinear SDR, methods more sophisticated than GSIR have been developed. For example, under the constant variance assumption that $\var \{\alpha(X)\mid f_0(X)\}$ is degenerate for all $\alpha(X) \in \lspf$, generalised sliced average variance estimation \citep[GSAVE;][]{lee13} detects $\alpha(X) \in \lspf$ using 
\begin{equation}\label{eq: gsave pop}
E \big[ \var\{\alpha(X)\mid Y\} - \var\{\alpha(X)\} \big]^2 = 0,
\end{equation}
and subsequently recovers a collection of transformations of $f_0(X)$. Similar to GSIR, these methods may also be reformulated algorithmically and generalised to become powerful dimensionality reduction methods for machine learning research. In the next section, we explore an example inspired by GSAVE.

\section{Second-order supervised Laplacian embedding}\label{sec: second-order sle}

\def\Gf{{\mathcal {G}}}
\def\Es{{\mathcal {E}}}

As seen in (\ref{eq: gsave pop}), GSAVE can be understood as finding a transformation of $X$ whose variance is maximally affected when conditioned on $Y$. Following the same spirit, we now recover the reduced predictor based on its higher-order alignment with $Y$; that is, the distance between two observations in the reduced space is maximally affected when conditioned on the corresponding responses. For a prespecified positive integer $d<p$, this amounts to maximising over $\psi(X) \in \Ldx d$
\begin{equation}\label{def: gsave equiv alt}
\var\big[E\{\|\psi(X) - \psi(X^*)\|^2 \mid Y,Y^*\}\big],
\end{equation}
where we replace the similarity measure $w(Y,Y^*)$ used in ordinary supervised Laplacian embedding with the ambient $(Y,Y^*)$ to fully exploit the information contained in this pair of responses.

Compared with (\ref{def: sle 1 pop}), which represents the average level of alignment between the reduced predictor and the response, (\ref{def: gsave equiv alt}) probes a subtler form of coherence between the two random objects. A potential issue with  (\ref{def: gsave equiv alt}) is that, in contrast to the reduced predictor derived from ordinary supervised Laplacian embedding, its maximiser may not satisfy the nesting property mentioned below Theorem \ref{thm: exhaustiveness}. For interpretability of dimensionality reduction, we modify (\ref{def: gsave equiv alt}) to sequentially extract the variables in the reduced predictor under their linear orthogonality constraints, thereby forcing the reduced predictor to satisfy the nesting property. 

To start with, we derive the first variable $\psi_{0,1} (X)$ by maximising over $\psi(X) \in \Ldx 1$
\begin{equation}\label{def: gsave equiv}
\Gf\{\psi(X)\} = \var\big[E[\{\psi(X) - \psi(X^*)\}^2 \mid Y,Y^*] \big].
\end{equation}
Then, in the $k$th step, we derive $\psi_{0,k} (X)$ by maximising $\Gf\{\psi(X)\}$ over $\psi(X) \in \Ldx 1$ subject to $E\{\psi(X) \psi_{0,\ell} (X)\} = 0$ for $\ell=1,\ldots, k-1$. To ensure the uniqueness of each $\psi_{0, k}(X)$, we impose the following regularity condition, similar to Condition (C2).

\vspace{0.2cm}

\noindent
(C5) The maximiser of $\Gf(\cdot)$ in (\ref{def: gsave equiv}) over $\psi(X) = g\circ f_0(X) \in \Ldfx 1$ is unique up to multiplication by $\pm 1$. The same holds for any $k > 1$ under the additional constraints $E[ \{ g\circ f_0(X)\} \{\phi_\ell \circ f_0(X)\}] =0$ for $\ell=1,\ldots, k-1$, where each $\phi_\ell \circ f_0(X)$ denotes the unique solution to the same maximisation problem with $k$ replaced by $\ell$. 

\vspace{0.2cm}

\noindent
Denote $(\psi_{0, 1} (X),\ldots, \psi_{0,d}(X))\trans$ by $\psi_0(X)$. By construction, we call any method that aims to deliver $\psi_0(X)$ second-order supervised Laplacian embedding. 



As with GSAVE, the effectiveness of second-order supervised Laplacian embedding is straightforward under the constant variance assumption mentioned above (\ref{eq: gsave pop}): if any $\alpha(X)$ in $\lspf$ induces a degenerate $E[\{\alpha(X) - \alpha(X^*)\}^2 \mid f_0(X),f_0(X^*)]$ and thus $\Gf\{\alpha(X)\}=0$, then the orthogonal complement of $\{\psi(X) \in \Ldx 1: \Gf\{\psi(X)\} = 0\}$ in $\Ldx 1$ must be a transformation of $f_0(X)$. Because we only need the weaker result that a small set of maximisers of $\Gf (\cdot)$ is a transformation of $f_0(X)$, we can relax the constant variance assumption to: 

\vspace{0.2cm}

\noindent
(C6) For any $\psi(X)\in \Ldx 1 \cap \lspf$, $\Gf\{\psi(X)\} < \Gf\{\psi_{0, k} (X)\}$ for all $k \geq 1$. 

\vspace{0.2cm}

\noindent
By the definition of $\Gf(\cdot)$, Condition (C6) roughly requires that any transformation of $X$ irrelevant to $Y$ cannot induce a metric of observations that varies drastically when conditioned on the responses.  
As shown in the Supplementary Material, this condition is implied by the completeness assumption mentioned in Section \ref{sec: sle theory}. Thus, it is fairly general and holds in particular for both homoscedastic regression and inverse regression models. Considering the impracticality of the constant variance assumption, this relaxation substantially enhances the applicability of second-order supervised Laplacian embedding.

\begin{theorem}\label{thm: gsave equiv}
Under Conditions (C5) and (C6), for any $k > 0$, $\psi_{0,k}(X)$ is unique up to multiplication by $\pm 1$, and must be a transformation of $f_0 (X)$. 
\end{theorem}



As with Theorem \ref{thm: exhaustiveness} for ordinary supervised Laplacian embedding, second-order supervised Laplacian embedding can be shown to fully recover $f_0(X)$ under appropriate regularity conditions. To this end, we speculate that Condition (C3) can be relaxed following the reasoning in \cite{lee13}, which demonstrates that GSAVE fully recovers $f_0(X)$ without the completeness assumption. The details are omitted. This relaxation can be regarded as a benefit of exploiting more sophisticated relationships between the reduced predictor and the response beyond their average alignment, in the same spirit as the belted neural network for nonlinear SDR \citep{tang2026belted} where a neural network is used to model the reduced data.
In practice, the additional flexibility of second-order supervised Laplacian embedding may provide better robustness to the choice of $\H_X$ compared with ordinary supervised Laplacian embedding. For example, if $(X,Y)$ exhibits the pattern in Figure \ref{fig: uneq wii}A in Section \ref{sec: sle theory}, then it is plausible that ordinary supervised Laplacian embedding will be effective only if $\H_X$ approximates an invertible transformation of $|X\comp 1|$ at the sample level, whereas second-order supervised Laplacian embedding additionally allows $\H_X$ to approximate an invertible transformation of $X\comp 1$ at the sample level (although $X\comp 1$ generates a larger $\sigma$-field than $\sigma_{Y|X}$ defined below (\ref{asm: sdr})). A potential price of second-order supervised Laplacian embedding is that its estimation may require a larger sample size due to the use of $E(\cdot \mid Y,Y^*)$. A thorough investigation of the bias--variance tradeoff between the two methods is deferred to future work. 


\section{Simulation studies}\label{sec: simulation}

We now use simulation models to evaluate the effectiveness of the proposed supervised Laplacian embedding method, in comparison to both GSIR and the earlier approach proposed in \cite{tai2022kernelized} (mentioned in the Introduction). To favour the latter two methods, we set the dimension of their reduced predictors to the true dimension of $f_0(X)$, known \textit{a priori}, and ensure the validity of this choice by carefully designing the simulation models to avoid the redundancy issue discussed in Section \ref{sec: sle theory}. In contrast, for the proposed supervised Laplacian embedding method, we use the BIC-type criterion proposed in Section \ref{sec: sle} to determine the dimension of its reduced predictor. The GMDD-based methods are omitted for two reasons: first, as discussed below (\ref{eq: gmdd}) in Section \ref{sec: sdr}, they are essentially special cases of our methods at the population level; second, they approximate $f_0 (X)$ by a deep neural network instead of an RKHS, which can affect their practical performance. Second-order supervised Laplacian embedding proposed in Section \ref{sec: second-order sle} is also omitted, as its implementation at the sample level requires further investigation.

For all three methods, we set $\H_X$ to the RKHS induced by the Gaussian kernel $K_X(X, X^*) = \exp \{-\|X-X^*\|^2/h_X^2\}$. For both the method in \cite{tai2022kernelized} and the proposed method, we set $w(Y,Y^*)$ to $\exp \{-(Y-Y^*)^2/h_Y^2\}$ if $Y$ is continuous and to $\ind (Y = Y^*)$ if $Y$ is discrete, both of which satisfy Condition (C1). For fairness of comparison, we use $w(Y,Y^*)$ also as the kernel function $K_Y(Y, Y^*)$ to construct $\H_Y$ in GSIR. As mentioned below Theorem \ref{thm: gsir reformulation}, this makes the proposed method nearly identical to GSIR when $Y$ is discrete, with negligible differences caused by the use of $\G_Y^\dagger$ in (\ref{eq: gsir sample}). The tuning parameters $\epsilon_n$ and $\epsilon_n^*$ in the approximations of $\G_X^\dagger$ and $\G_Y^\dagger$ are set uniformly at $0.01$. Both $h_X$ in $K_X(X, X^*)$ and $h_Y$ in $w(Y,Y^*)$ (for continuous $Y$) are tuned for each individual method, using five-fold cross-validation where the goodness of fit on the test sets is measured by the sample distance correlation between the reduced predictor and $Y$.

We generate data from each of the following four models, where $\varepsilon\comp 1$ and $\varepsilon\comp 2$ are random errors distributed as $N(0, (1/3)^2)$, $\Sigma_0$ denotes the $p\times p$ matrix that equals $\I_p$ except for its $(i,i+1)$th and $(i+1,i)$th entries being $0.5$ uniformly for $i=1,\ldots,p-1$, and $U_p(-2,2)$ denotes the distribution in $\real^p$ under which all variables are independent and uniformly distributed on $(-2, 2)$.

\vspace{0.2cm}

\noindent
Model 1: $X \sim N(\zero, \Sigma_0)$ and $Y\mid X \sim \mathrm{Bernoulli}\{\mathrm{expit}(3 X\comp 1 - 3 X\comp 2)\}$.

\vspace{0.2cm}

\noindent
Model 2: $X \sim N(\zero, \I_p)$ and $Y = |X\comp 1 | + |X\comp 2 + 1 | + \varepsilon\comp 1$. 

\vspace{0.2cm}

\noindent
Model 3: $X \sim U_p(-2, 2)$ and $Y = {\mathrm {sign}}(X\comp 1 + \varepsilon\comp 1) + {\mathrm {sign}}(X\comp 2 + \varepsilon\comp 2)$. 

\vspace{0.2cm}

\noindent
Model 4: $X \sim U_p(-2, 2)$ and $Y = \{X\comp 1\}^2 + {\mathrm {sign}}(\varepsilon\comp 1) \{X\comp 2\}^2 + \varepsilon\comp 2$. 

\vspace{0.2cm}

\noindent
In Models $1$ and $3$, $Y$ is discrete. In Models $2$ and $4$, $Y$ is continuous. The desired reduced predictor $f_0(X)$ is univariate in Models $1$ and $2$ and bivariate in Models $3$ and $4$. Because a larger-dimensional $f_0(X)$ is usually more difficult to estimate, we set $(n,p) = (200, 10)$ in Models $1$ and $2$ and increase $n$ to $400$ in Models $3$ and $4$. For each model, we generate $500$ independent replicates  and apply the three aforementioned dimensionality reduction methods. The deviation of a reduced predictor $\widehat f_0(X)$ from the desired $f_0(X)$ is measured by their sample distance correlation calculated on an independent test set of size $5000$. This measure ranges from $0$ to $1$, with a larger value indicating a smaller deviation of $\widehat f_0(X)$ from $f_0(X)$. Based on $500$ repetitions, the results are summarised in Table \ref{table: sim est}, where TSLE denotes the supervised Laplacian embedding method proposed in \cite{tai2022kernelized}.


\begin{table}[h]
\begin{center}
\centering
\caption{Effectiveness of dimensionality reduction methods: ``TSLE'' denotes the supervised Laplacian embedding method in \cite{tai2022kernelized} and ``Proposed'' denotes the new method proposed in Section \ref{sec: sle}. In each cell of Columns $2$-$4$, $a(b)$ is the sample mean (sample standard deviation) of the distance correlation between the corresponding reduced predictor and $f_0(X)$, evaluated based on $500$ replications. In each cell of Column $5$, the number is the percentage of times that the proposed BIC-type criterion correctly specifies the dimension of the reduced predictor for the proposed method, evaluated based on $2000$ replications.}
\begin{tabular}[l]{ccccc}
\hline \hline
\empty    &   TSLE          & GSIR          & Proposed  & BIC \\ 
\hline
Model $1$ &  $0.26(0.09)$   & $0.90(0.02)$  & $0.90(0.02)$ & $99$ \\ 
Model $2$ &  $0.12(0.04)$   & $0.85(0.02)$  & $0.83(0.02)$ & $95$ \\     
Model $3$ &  $0.40(0.08)$   & $0.77(0.01)$  & $0.77(0.01)$ & $100$ \\  
Model $4$ &  $0.11(0.02)$   & $0.81(0.03)$  & $0.84(0.02)$ & $99$ \\
\hline\hline
\end{tabular}
\label{table: sim est}
\end{center}
\end{table}

\begin{figure}[h]
\hspace*{-0.75cm} 
     \centering
            \includegraphics[page=1, width=0.36\textwidth]{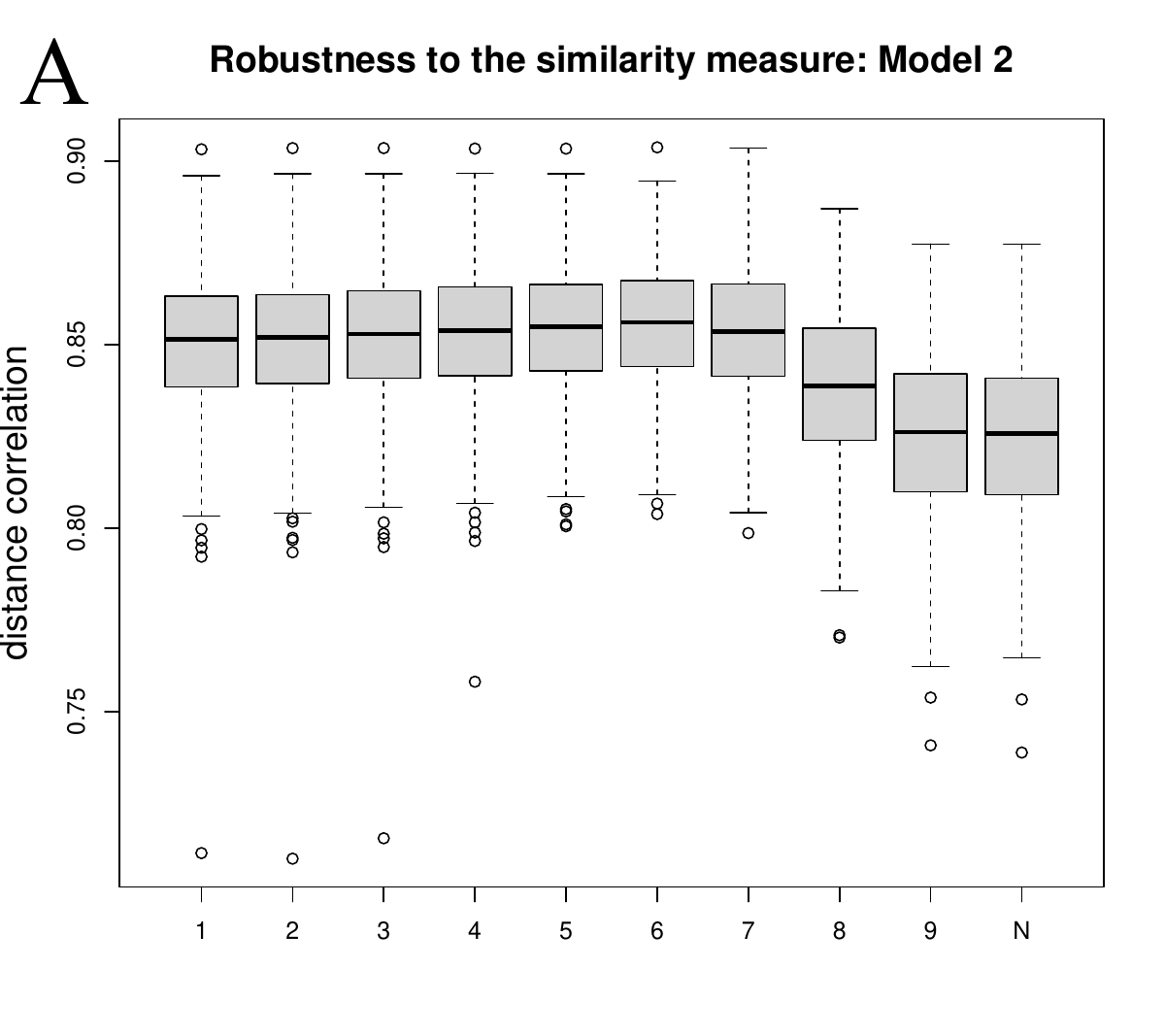}
            \includegraphics[page=2, width=0.36\textwidth]{tune_gsir.pdf}
    \caption{Performance of GSIR in Model $2$ (Panel A) and Model $4$ (Panel B) with the tuning parameter $\epsilon_n^*$ set to $10^{i-7}$ for $i=1,\ldots, 10$. The numbers on the $x$-axis indicate the values of $i$, while ``$N$'' on the $x$-axis corresponds to the proposed method. Each boxplot summarises the performance of GSIR (or the proposed method), measured by the sample distance correlation between the reduced predictor and $f_0(X)$ based on a test set of size $5000$.}
    \label{fig: tune gsir}
\end{figure}

From Table \ref{table: sim est}, the supervised Laplacian embedding method proposed in \cite{tai2022kernelized} fails to deliver the desired reduced predictor in all models, which agrees with the theoretical discussion in the Introduction. By contrast, the proposed method is effective in all models. As the rightmost column in Table \ref{table: sim est} shows, this effectiveness is partially due to the consistency of the proposed BIC-type criterion in determining the dimension of the reduced predictor. 
Interestingly, while GSIR and the proposed method are theoretically anticipated to be nearly identical in Models $1$ and $3$, where $Y$ is discrete, their performance is also indistinguishable in Models $2$ and $4$, where $Y$ is continuous. This validates the robustness of GSIR to the tuning parameter $\epsilon_n^*$ mentioned in Section \ref{sec: sdr}, as the proposed method approximates an extreme case of GSIR where an excessively large $\epsilon_n^*$ essentially replaces $\G_Y (\G_Y + \epsilon_n^* \I_n)\inv$ with $\G_Y / \epsilon_n^*$. To investigate this more closely, we assess GSIR in Models $2$ and $4$ with $\epsilon_n^*$ set to $10^{i-7}$ for $i=1,\ldots, 10$. The results, summarised in Figure \ref{fig: tune gsir}, suggest that GSIR performs nearly the same as the proposed method regardless of the choice of $\epsilon_n^*$. 

\section{A real data example}\label{sec: real data}

\begin{figure}[h]
\hspace*{-0.75cm} 
     \centering
            \includegraphics[page=1, width=0.32\textwidth]{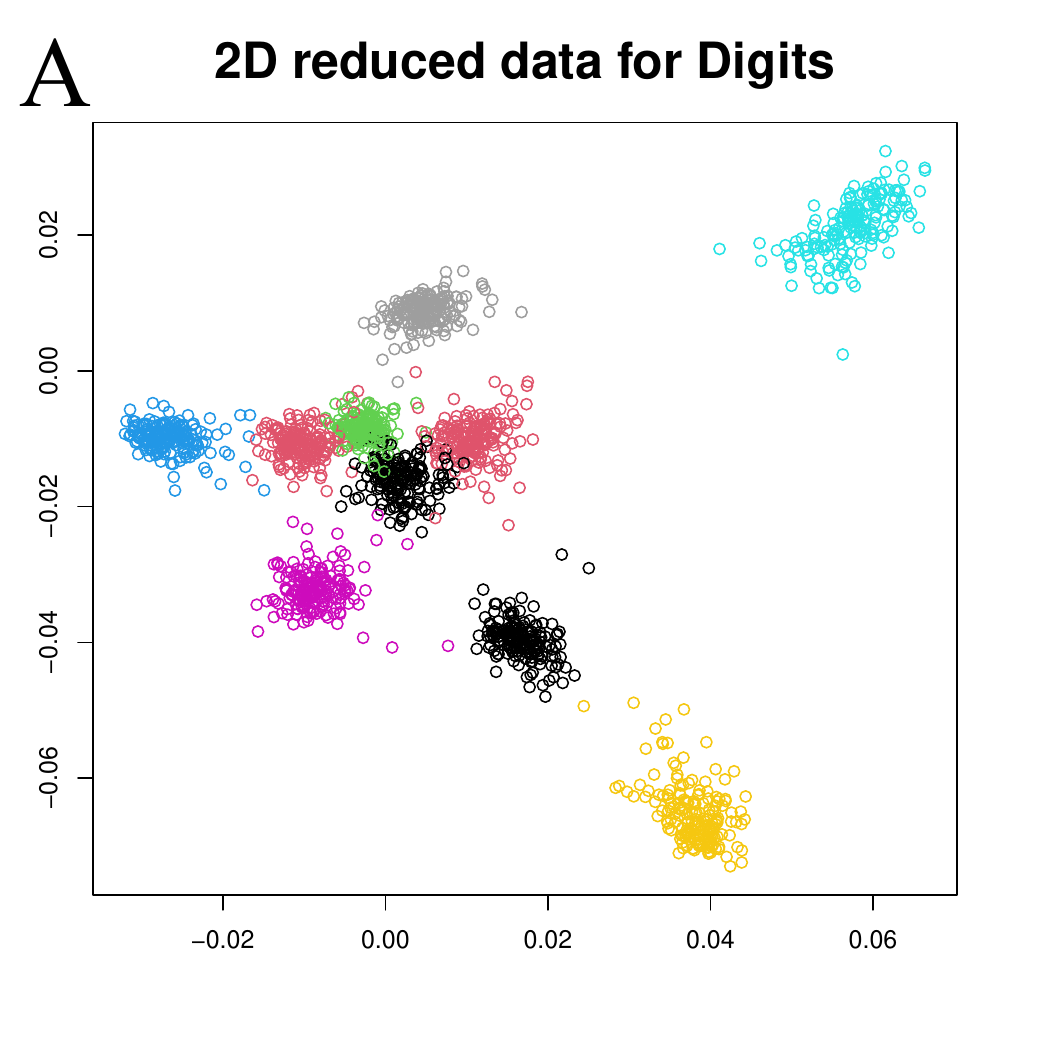}
            \includegraphics[page=2, width=0.32\textwidth]{digits.pdf}
            \includegraphics[page=3, width=0.32\textwidth]{digits.pdf}
\caption{Panel A shows the scatterplot of the reduced predictor for the Digits dataset obtained from the proposed supervised Laplacian embedding method, with $d=2$,  x- and y-axes corresponding to the first and second variables, and colors indicating different digits. Panel B shows the curve of the test misclassification rate of linear discriminant analysis for the proposed reduced data evaluated via $10$-fold cross-validation, where the dimension of the reduced predictor varies from $1$ to $10$. Panel C shows the same scatterplot as in Panel A but with $d=3$, x-, y-, and z-axes corresponding to the first, third, and second variables, and only the observations for digits $1$, $2$, $8$, and $9$ are displayed. }
    \label{fig: digits}
\end{figure}

The Digits dataset, available at 
\url{https://scikit-learn.org/stable/modules/generated/sklearn.datasets.load_digits.html}, contains $1797$ images of handwritten digits ranging from zero to nine, with a roughly equal number of images per digit. Each image is $8\times 8$ grayscale, yielding a $64$-dimensional $X$. Using a Gaussian $\H_X$ and setting $w(Y,Y^*) = \ind (Y=Y^*)$, we apply the proposed supervised Laplacian embedding method to this dataset. The bandwidth $h_X$ for $\H_X$ is set at $75$, tuned via $10$-fold cross-validation, where linear discriminant analysis is applied to the reduced data in the test sets with goodness of fit measured by the misclassification rate. Assuming that the reduced predictor is at least bivariate, we plot its first two variables in Figure \ref{fig: digits}A, with each digit class shown in a different color. From Figure \ref{fig: digits}A, the proposed reduced predictor clearly distinguishes most of the digits. To verify that this is not an artifact of overfitting, we conduct the $10$-fold cross-validation described above, and find that the misclassification rate of linear discriminant analysis on the reduced data is $0.12$. Figure \ref{fig: digits}A also shows that the distribution of the bivariate reduced predictor within each digit class is approximately normal with a homogeneous covariance matrix, supporting the use of linear discriminant analysis for subsequent classification.

For the dimension of the reduced predictor, the proposed BIC-type criterion suggests $d=9$. To evaluate this choice, we apply linear discriminant analysis to the proposed reduced data where the dimension of the reduced predictor is set to each of $1,\ldots, 10$, and approximate the corresponding test misclassification rate using the aforementioned $10$-fold cross-validation. The results plotted in Figure \ref{fig: digits}B show that the classification accuracy improves monotonically as the dimension of the reduced predictor increases from $1$ to $9$ and remains unchanged from $9$ to $10$, suggesting the optimality of taking $d=9$. Under the concern of model parsimony, the modest improvement from $3$ to $9$ suggests that taking $d=3$ is also a reasonable choice, which lowers the test misclassification rate to $0.06$ compared with $0.12$ for $d=2$ mentioned above. To illustrate the gain from increasing $d$ from $2$ to $3$, we pick the subsample associated with the four digits $\{1,2,8,9\}$, which are located near the centre of the scatterplot of the bivariate reduced predictor, and depict the three-dimensional reduced predictor for this subsample in Figure \ref{fig: digits}C. Comparing the two scatterplots in Figure \ref{fig: digits}, the third variable in the reduced predictor clearly helps to separate these four digit classes.


For reference, \cite{tai2022kernelized} also analysed this dataset, setting $h_X$ to $(0.0015)^{-1/2}$ for a Gaussian $\H_X$, the balance parameter $\lambda$ (mentioned in the Introduction) to $0.3$, and the dimension of the reduced predictor to two. The resulting reduced predictor, visualized in their Figure 8, exhibits a pattern similar to ours in Figure \ref{fig: digits}. However, based on the same $10$-fold cross-validation described above, their misclassification rate is as high as $0.51$, suggesting that the proposed supervised Laplacian embedding method can be particularly effective in practice.

\bibliographystyle{abbrv}
\bibliography{ref}

\end{document}